\input amstex
\documentstyle{amsppt}
\vsize=7.0in \hsize 6.8 in
\voffset=2cm \loadbold \topmatter
\title Rigidity of non-negatively curved metrics on open five-dimensional manifolds\endtitle
\author Valery Marenich \endauthor
\thanks
\endthanks
\address H\"{o}gskolan i Kalmar, 391~82, Kalmar, Sweden\endaddress \email valery.marenich\@ hik.se \endemail

\abstract
As the first step in the direction of the Hopf conjecture on the non-existence of metrics with positive sectional curvature on
$S^2 \times S^2$ the authors of [GT] suggested the following (Weak Hopf) conjecture (on the rigidity of non-negatively curved
metrics on $S^2 \times R^3$): "The boundary $S^2\times S^2$ of the $S^2 \times B^3\subset S^2 \times R^3$ with an arbitrary
complete metric of non-negative sectional curvature contains a point where a curvature of $S^2 \times S^2$ vanish". In this
note we verify this.
\endabstract
\subjclass   53C20, 53C21. Supported by the Swedish Science Consul ({Vetenskapr{\aa}det}) and the Faculty of Natural Sciences
of the Hogskolan i Kalmar, (Sweden). Submitted August 12, 2004; revised November 29, 2004
\endsubjclass
 \dedicatory  \enddedicatory
\endtopmatter

\document
\baselineskip 12pt

\head \bf More "flats" in $M^{5}$ \endhead

Let $(M^{n},g)$ be a complete open Riemannian manifold of non-negative sectional curvature. Remind that as follows from [CG]
and [P] an arbitrary complete open manifold $M^n$ of non-negative sectional curvature contains a closed absolutely convex and
totally geodesic submanifold $\Sigma$ (called a soul) such that the projection $\pi: M\to \Sigma$  of $M$ onto $\Sigma$ along
geodesics normal to $\Sigma$ is well-defined and is a Riemannian submersion.\footnote{of class $C^{1,1}$. Some additional
efforts should be made to verify that $\pi$ is of the same class of smoothness as $M$ in order to have O'Neill's fundamental
equations for Riemannian submersions, see Appendix A below for this and some other results.} The (vertical) fibers
$W_P=\pi^{-1}(P), P\in \Sigma$ of $\pi$ define a metric foliation in $M$ and two distributions: a vertical ${\Cal V}$
distribution of subspaces tangent to fibers and a horizontal distribution ${\Cal H}$ of subspaces normal to ${\Cal V}$. For an
arbitrary point $P$ on $\Sigma$, an arbitrary geodesic $\gamma(t)$ on $\Sigma$ and arbitrary vector field $V(t)$ which is
parallel along $\gamma$ and normal to $\Sigma$ the following
$$
\Pi(t,s)=exp_{\gamma(t)} sV(t) \tag 1
$$
are totally geodesic surfaces in $M^n$ of zero curvature, i.e., flats.

Sometimes, these are the only directions of zero curvature in open $M^n$ (e.g., when $M^4$ is the tangent bundle to the
two-dimensional sphere with the Cheeger-Gromoll metric, see [M2]). The objective of this note is to verify the (Weak Hopf)
conjecture from [GT] and to point to more directions of zero curvature in our particular case of a five-dimensional $M$. The
following statement is true.

\proclaim{Theorem~A} There does not exist a metric of nonnegative curvature on $M^5$ for which the boundary of a small metric
tube about the soul has positive curvature in the induced metric.
\endproclaim

Clearly, the only difficult case in the Theorem~A is of two-dimensional $\Sigma$ diffeomorphic to $S^2$. For other cases of
$codim(\Sigma)=1,2$ or $4$, or two-dimensional and non-orientable $\Sigma$ or torus might be easily treated or by going to the
oriented covering, or by applying "the straight line splitting off" theorem by Toponogov. Note also that unlike [GT] we are not
assuming that the normal bundle of the soul is topologically trivial.\footnote{Which is not really a strong restriction since
there are only two non-homotopic vector bundles over $S^2$: trivial and non-trivial for which the corresponding unit-sphere
bundle is a ruled surface - the only non-trivial $S^2$-bundle over $S^2$. These bundles correspond to elements of the
$\pi_1(SO(3))=Z_2$ and both admit a non-zero section. Thus, our main technical result, Theorems~2a and 2b below may be
considered as yet another splitting result: the local existence of the parallel sections when the curvature of $M$ is
non-negative.}

The proof of Theorem~A is based on the consideration of some family of holonomy operators in $M$.\footnote{and is a further
development of our "prism" construction from [M1,3].} More precisely, we consider a disk $\Omega$ in $\Sigma$ bounded by a
curve $\omega$, construct a smooth homotopy $\omega_x$ of this curve to a point and consider the family of parallel transports
$I_{\omega_x}$ along $\omega_x$ acting on vectors normal to $\Sigma$. Our construction heavily depends on $dim(\Sigma)=2$ and
$codim(\Sigma)=3$ conditions which makes its generalizations to higher dimensions difficult.

The proof of the Theorem~A is given in the section~5 after the construction of the family of holonomy operators in the
section~1, curvature calculations in the section~2, constructions of the local (and the global) parallel sections in the
section~3 (and 4 correspondingly).


\head \bf 1. The Holonomy and the O'Neill's $A$-tensor \endhead

Let $\Omega$ be a disk in two-dimensional sphere $\Sigma$ bounded by a closed curve $\omega$. According to the construction
given below (see subsection 1.3) $\Omega=\{\cup \omega_x(y)| 0<x\leq 1\}\cup \{O\}$, where $O$ is some interior point
("center"), the boundary curve $\omega=\partial \Omega$ equals $\omega_1$, $\omega_x(y),0\leq y\leq 2\pi$ is the family of
closed curves such that $\{x,y | 0<x\leq 1, 0\leq y<2\pi\}$ are ("polar-type") "coordinates" in $\Omega\backslash \{O\}$. The
point with "coordinates" $\{x,y\}$ we denote by $P(x,y)$, and do not assume that the correspondence $(x,y)\to P(x,y)$ is
one-to-one. We also assume that the parameter $y$ on $\omega_x$ is proportional to the arc-length. Let $X(x,y)$ and $Y(x,y)$ be
an orthonormal base of $T_{P(x,y)}\Sigma$ with positive orientation such that unit $Y(x,y)$ has the same directions as
$\dot\omega_x(y)=\partial P(x,y)/\partial y$.

Fix some positive $s_0$ smaller than a focal radius of $\Sigma$ in $M$. For some $s<s_0$ denote by $N\Sigma(s)$ the boundary of
an $s$-neighborhood of $\Sigma$. Due to our choice it is a smooth manifold. It consists of all points $Q(P,V)=exp_P(sV)$, where
$P$ is a point on $\Sigma$ and $V$ is a unit vector normal to $\Sigma$ at $P$. A unit normal $V(Q)$ to hyper-surface
$N\Sigma(s)$ at $Q(P,V)$ is the parallel translation of $V$ from $P$ to $Q$ along a vertical geodesic $exp_P(s'V), 0\leq s'\leq
s$. By $\bar X(x,y,s,V)$ and $\bar Y(x,y,s,V)$ (or simple $\bar X(x,y)$ and $\bar Y(x,y)$ if there is no confusion) we denote
horizontal lifts of $X(x,y)$ and $Y(x,y)$ from $P=P(x,y)$ to $Q(P,V)$.

By a vertical lift of a point $P\in\Sigma$ in direction $V\in\nu_P\Sigma$ we mean a point $Q=exp_P(sV)$ for some $s>0$.
Correspondingly, $\omega(y,V(y))=exp_{\omega(y)}(sV(y))$ is said to be a vertical lift of $\omega(y)$ along some vertical
vector field $V(y)$ along $\omega$. Due to (1) when $V(y)$ is a parallel vertical vector field along $\omega$ its vertical lift
$\omega(y,V(y))$ is a horizontal curve (i.e., its speed is a horizontal vector everywhere). In this case we say, as usual, that
$\omega(y,V(y))$ is a horizontal lift of $\omega$ (see [O'N]). The map $\pi: \omega(y,V(y))\to\omega(y)$ decrease the distance
(i.e., is "short") and is an isometry iff $\omega(y,V(y))$ is a horizontal lift of $\omega$.

The total vertical lift of $\omega_x$, i.e., the sub-manifold $\Psi_x(y,V)=exp_{\omega_x(y)}sV$ for all unit $V\in
\nu_{\omega_x(y)}\Sigma$ is a collection over $\omega_x$ of all vertical $s$-spheres. If $X(x,y)$ is a unit vector in
$T_{\omega_x(y)}\Sigma$ normal to $\dot\omega_x(y)$ then its parallel transport $\bar X(Q)$ from $\omega_x(y)$ to
$Q=Q(\omega_x(y), V)$ (along vertical geodesic) is a normal to $\Psi_x$, so that $\{V(Q), X(Q)\}$ is an orthonormal base of the
normal subspace to $\Psi_x$ at $Q$.

\subhead 1.1 O'Neill's fundamental equations \endsubhead

Remind, that according to the fundamental O'Neill's formula (see [O'N])\footnote{and also [M1-3] for an exposition adapted to
our case.}
$$
(R(X,V)Y,W)=((\nabla_X T)_V W,Y) + ((\nabla_V A)_X Y,W) - (T_V X,T_Y W) + (A_X V,A_W Y) \tag 2
$$
where $X,Y$ are horizontal vectors (i.e., belonging to ${\Cal H}$), $V,W$ are vertical vectors (i.e., belonging to ${\Cal V}$)
and $T$ and $A$ are O'Neills fundamental tensors defined as follows
$$
T_E F={\Cal V}(\nabla_{{\Cal V}(E)}({\Cal H}(F)))+{\Cal H}(\nabla_{{\Cal V}(E)}({\Cal V}(F))) \qquad \text{ and }\qquad A_E
F={\Cal V}(\nabla_{{\Cal H}(E)}({\Cal H}(F)))+{\Cal H}(\nabla_{{\Cal H}(E)}({\Cal V}(F))).
$$
Here tensor $T$ is the second fundamental form of vertical fibers, while $A$ measures non-integrability of the horizontal
distribution. Therefore,
$$
(R(X,W)W,X)=((\nabla_X T)_V W,X) - \|(T_WX\|^2 + \|A_X W\|^2 \tag 3
$$
because, as easy to verify,
$$
((\nabla_W A)_X X,W)=0
$$
due to the fact that $A$ is anti-symmetric and horizontal, see again [O'N]. Vanishing of the curvature term in (3) will imply
below Theorem~A. Another fundamental formula by O'Neill:
$$
(R(X,Y)Y,X)(P)=(R(X,Y)Y,X)(Q(P,V)) + 3\|A_{X}Y\|^2(Q(P,V)). \tag 4
$$

\subhead 1.2 Prism construction \endsubhead

From [M1-3] we have the following.

\proclaim{Lemma~1}
$$
\|A_{X}Y\|^2(Q(P,V))= {{s^2}\over{4}} \|R(X,Y)V\|^2(P).
$$
\endproclaim

The sketch of the proof of Lemma~1 is (see [M1-3] for calculations): take a small triangle $\triangle PP_1P_2$ with sides
parallel to $X$ and $Y$, translate parallel $V$ along these sides to vectors $V_1$ and $V_2$ at $P_1$ and $P_2$ correspondingly
and lift-up the vertices of the triangle in obtained directions: $\triangle(s')=\triangle P(s)P_1(s)P_2(s)$ (we have a "prism")
where $P(s)=exp_{P}sV$, $P_1(s)=exp_{P_1}sV_1$, $P_2(s)=exp_{P_2}sV_2$. From (1) it follows that the angle $\angle P(s)$ and
sides of $\triangle(s)$ have zero first and second derivatives. Hence, the second derivative of the length of the third side
$P_1(s)P_2(s)$ is proportional to the second derivative of the curvature of $M$ in two-dimensional direction $\{X,Y\}$. The same
second derivative of the length of the third side $P_1(s)P_2(s)$ can be computed in a different way: by comparing $V_2$ with the
parallel translation $V'_2$ of $V_1$ from $P_1$ to $P_2$ along $P_1P_2$. By Ambrose-Singer theorem $V_2-V'_2$ translated from
$P_2$ to $P$ equals $R(X,Y)V$ times the area of the triangle $\triangle PP_1P_2$ up to higher order terms. Then the second
variation formula due to (1) implies the claim of the Lemma~1.

Before going further remark, that $A_XY(Q)$ does not depend on the particular choice of the orthonormal base $X,Y$ with a
positive orientation of a horizontal subspace ${\Cal H}_Q$. Indeed, due to $A_HH\equiv 0$ for another orthonormal base with a
positive orientation $\tilde X=cos(\alpha)X+sin(\alpha)Y, \tilde Y=-sin(\alpha)X+cos(\alpha)Y$ we have
$$
A_{\tilde X}\tilde Y=(cos^2(\alpha)+sin^2(\alpha)) A_X Y =A_X Y.
$$
Therefore, in what follows we denote $A_X Y(Q(P,V))$ simply $A(Q)$ for $Q=Q(P,V)$.

Vanishing of $A$ implies
$$
A_XW\equiv 0 \tag 5
$$
for all horizontal $X$ and vertical $W$, i.e., that the vertical subspace is parallel in horizontal direction. Indeed,
$$
A_XW=(A_XW,X)X+(A_XW,Y)Y=(A_XW,Y)Y=({\Cal H}(\nabla_{{\Cal H}(X)}({\Cal W}(F))),Y)=-(W,A_XY).
$$

Because $A(Q)$ is orthogonal to the normal $V(Q)$ of $N\Sigma(s)$ and vertical, it defines a vector field tangent to the
vertical two-dimensional sphere $S^2(P)(s)=N\Sigma(s) \cap W_P$. Therefore, $A(Q)$ vanish at some $Q^*=Q(P,V^*(P,s))$ for every
$P$. Note that from the Lemma~1 we deduce:

\proclaim{Lemma~2} For a given $P$ the vector $V^*(P,s)$ does not depend on $s$ and satisfies $R(X,Y)V^*=0$. For a fixed $s$ the
set of all $Q=Q(P,V,s)=exp_P sV$ such that $A(Q)=0$ is in one-to-one correspondence with the set of $V\in \nu_P\Sigma$ such that
$R(X,Y)V=0$.
\endproclaim

Now we employ $codim \Sigma=3$.

As we saw $A(Q)$ is proportional to the generator $R(X,Y)V$ of the holonomy group of the normal bundle $\nu\Sigma$, and
therefore $A$ or is identically zero on a given vertical sphere $S^2(P)$, or vanish for two opposite to each other normals
$V^*_1$ and $V^*_2=-V^*_1$, which by the Lemma~2 does not depend on the radius $s$ of vertical spheres; while parallel
translations of the space $\nu_P\Sigma$ normal to $\Sigma$ around small closed contours around $P$ in positive direction are
rotations about $V^*_1$ in positive direction with a speed equals to the area bounded by the contour times $|(R(X,Y)W,U)|$
where $W,U$ from $\nu_P\Sigma$ are orthonormal and orthogonal to $V^*_1$. If we denote by $Hol(P)$ a rotation of $\nu_P\Sigma$
about an axis $V^*_1$ in positive direction and speed $|R(X,Y)W,U|$ (a density of the holonomy operator according to
Ambrose-Singer theorem) we will have a continuous map $Hol:\Sigma\to SO(3)$, which we call an infinitesimal holonomy map - a
nice geometric representation for the holonomy of the normal bundle $\nu\Sigma$.

\subhead 1.3 Construction of the homotopy $\omega_x$ \endsubhead

In [M1] (see also [M3,4]) we proved that if the holonomy of the normal bundle $\nu\Sigma$ of the simply connected soul in an
open manifold $M^n$ of non-negative curvature is trivial then the manifold $M^n$ is isometric to the direct product. In this
case the Theorem~A is obviously true. Thus, we may assume that at some point $O\in\Omega$\footnote{change $\Omega$ to
$\Sigma\backslash \Omega$ if necessary} the generator of the holonomy operator is not zero, i.e., $R(X,Y)V(O)\not\equiv 0$ so
that the vector $V^*(O)$ as above is uniquely defined. Our construction we start with some initial homotopy $\omega_x, 0\leq
x\leq 1$ of $\omega=\partial \Omega$ to a point $O$, i.e., such that $\Omega=\{\cup \omega_x(y)| 0<x\leq 1\}\cup \{O\}$; and
then will change it if necessary.

Consider the parallel translation $I_x$ of $\nu_{P(x)}\Sigma$, where $P(x)=\omega_x(0)$, into itself along $\omega_x$ - we call
it the holonomy along $\omega_x$. Due to our choice of $O$ it is not the identity map for small $x$, and because $codim
\Sigma=3$ this holonomy is a rotation about some uniquely defined axis generated by a vector $V(x)\in \nu_{P(x)}\Sigma$ such
that $V(x)\to V^*(O)$ as $x\to 0$. For definiteness we choose $\omega_x$ equals a circle of radius $x$ around $P$ for small $x$.
Then $I_x$ depends smoothly on $x$, and because $V(x)$ is uniquely defined - it also depends smoothly on $x$ for sufficiently
small $x$. Then the image $V(x,y)$ of the vector $V(x)$ under the parallel transport $I_x(y)$ along $\omega_x$ from $P(x)$ to
$\omega_x(y)$ is also a smooth vector field. This will imply that the surfaces $\Omega_s$ which we will construct below as
vertical lifts of $\Omega$ in the direction of the vector field $V$ will be smooth. Note that it always holds
$$
\nabla_Y V(x,y)\equiv 0,
$$
and it is not difficult to see that all first covariant derivatives of $V(x,y)$ actually vanish at $O$.

Consider how $I_x$ varies for bigger $x$. If for all $0<x\leq 1$ it is a rotation on non-zero angle about some uniquely defined
vector $V(x)$ we have our homotopy defined. Otherwise for some $x\nearrow x^*$ the family of holonomies $I_x$ converges (in a
natural sense) to $I_{x^*}$ which is the identity map. In other words, if $H:(0,1]\to SO(3)$ is the action of $I_x$ on
$\nu_O\Sigma$ as follows:
$$
H(x)(V)=J^{-1}(x)\circ I_x\circ J(x)(V),
$$
where $J(x)$ is a parallel translation from $O$ to $P(x)$, then $H(x)\in SO(\nu_P\Sigma)=SO(3)$ and $H(x^*)=id$. Having this
trouble, i.e., $H(x^*)=id$ we may try to "take off" the curve $H(x), x^*-\delta<x<x^*+\delta$ of orthogonal transformations from
an identity point $id$ in $SO(3)$ by varying "the curve of curves" - the family $\omega_x$, i.e., taking some variation
$\tilde\omega_{x,\epsilon}$, where $\tilde\omega_{x,0}=\omega_x$ such that the new holonomy curve $H(\epsilon, x)
=J^{-1}(x)\circ I_{x,\epsilon}\circ J(x)$ where $I_{x,\epsilon}$ is a parallel translation along $\tilde\omega_{x,\epsilon}$,
does not go through $id$ in $SO(3)$.

To simplify forthcoming computations we consider variations given by
$$
\tilde\omega_x(\epsilon,y)=exp_{\omega_{x}(y)} (\epsilon\phi_x(y)X(x,y)), \qquad x^*-\delta<x<x^*+\delta \tag 6
$$
where $X(x,y)$, as above, is a unit vector normal to $Y(x,y)$; and a smooth function $\phi_x(y)$ satisfying restrictions:
$$
\phi_x(y)\equiv 0 \quad\text{ for } \quad x < x^*-\delta, x^*+\delta < x.
$$
The varied family of holonomies along $\tilde\omega_x(\epsilon,y)$ defines the map $H(x,\epsilon)$ on two variables (and
depending on the "profile function" $\phi$) into three-dimensional $SO(3)$ as follows:
$$
H(x,\epsilon)(V)=J^{-1}(x,\epsilon)\circ I_x(\epsilon)\circ J(x,\epsilon)(V),
$$
where $J(x,\epsilon)=\tilde J(\epsilon)\circ J(x)$ and $\tilde J(\epsilon)$ is a parallel translation from $P(x)$ to
$P(x,\epsilon)=\tilde\omega_x(\epsilon,0)$ along $\tilde\omega_x(\epsilon,0)$, and $I_x(\epsilon)$ is a parallel translation
along $\tilde\omega_x(\epsilon,y)$ of vectors from $\nu_{P(x,\epsilon)}$. If partial derivatives of $H(x,\epsilon)$ on $x$ and
on $\epsilon$, i.e., two vectors $\partial H(x,\epsilon) /\partial x$ and $\partial H(x,\epsilon) /\partial \epsilon$ are
linearly independent at $(x^*,0)$ for a given $\phi$ then, obviously, there exists a variation $\tilde\omega_{x,\epsilon}$ for
which the curve $H(x,\epsilon)$ does not go through the point $id\in SO(3)$ for sufficiently small $\epsilon$. By Ambrose-Singer
theorem the action of the derivative $\partial H(x,\epsilon) /\partial \epsilon$ on a vector $W$ from $\nu_{O}\Sigma$ is
$$
(\partial H(x,\epsilon) /\partial \epsilon)(W) = J^{-1}(x^*)(\int\limits_{\omega_x} \phi_x(y)I^{-1}_x(y) R(X,\partial/\partial
y)W_x(y) dy)=
$$
$$
J^{-1}(x^*)({{L(x)}\over{2\pi}}\int\limits_{\omega_x} \phi_x(y)I^{-1}_x(y) R(X,Y)W_x(y) dy), \tag 7
$$
where $L(x)$ is the length of $\omega_x$, $I_x(y)$ denotes the parallel translation along $\omega_x$ from $\omega_x(0)$ to
$\omega_x(y)$ and $W_x(y)=I_x(y)(J(x^*)W)$. In particular for $\partial H(x,\epsilon) /\partial x$ it holds
$$
(\partial H(x,0) /\partial x)(W) = J^{-1}(x^*)({{L(x)}\over{2\pi}}\int\limits_{\omega_x} \psi_x(y) I^{-1}_x(y) R(X,Y)W_x(y) dy)
\tag 8
$$
where $\psi_x(y)=(X,\partial /\partial x)$.

We consider $R(X,Y)$ as an operator $R(X,Y): V\to R(X,Y)V$ from the Lie-algebra of $SO(3)$ which generates the holonomy group
and obtain some conditions on these generators in the case when both the first and the second differentials of $H(x,\epsilon)$
are degenerated at $(x^*,0)$ which do not allow to "take off" the curve $H(x)$ of the point $id$ in $SO(3)$. For short we
denote below $d_x H(x)= \partial H(x,\epsilon) /\partial x$ and $\delta_\epsilon H=\partial H(x,\epsilon) /\partial \epsilon$
and consider two possibilities for the vector $R=d_x H(x^*,0)$:\footnote{The same consideration might be done in coordinates:
choose an orthonormal base $\{E_i, i=1,2,3\}$ in $\nu_O\Sigma$ and the corresponding standard basis of $so(3)$ consisting of
three generators $R_{12},R_{13},R_{23}$ of $so(3)$ which are unit tangents to rotations in $\nu_O\Sigma$ with axes
$E_3,E_2,E_1$ correspondingly. Define vector fields $E_i(x,y)=I_{x}(y)J(x)E_i$ along $\omega_{x}$. Because $H(x^*,0)=id$ these
are (continuous) parallel vector fields along $\omega_{x^*}$. If $H(x,\epsilon)$ is given by the matrix $(H_{ij}(x,\epsilon);
i,j=1,2,3)$ then its derivatives have the following components:
$$
{{\partial H_{ij}(x^*,\epsilon)}\over{\partial \epsilon}}_{|\epsilon=0} ={{L(x)}\over{2\pi}} \int\limits_{\omega_{x^*}}
\phi_x(y) (R(X,Y)E_i(y),E_j(y)) dy. \tag 9
$$
where $E_i(y)=E_i(x^*,y)$ and
$$
{{\partial H_{ij}(x,0)}\over{\partial x}} = {{L(x)}\over{2\pi}}\int\limits_{\omega_x}\psi_x(y)(R(X,Y)E_i(y),E_j(y))dy. \tag 10
$$
These components are coordinates of derivatives of $H(x,\epsilon)$ in the basis $R_{12},R_{13},R_{23}$ of $so(3)$:
$$
\partial H(x,\epsilon) /\partial \epsilon = \sum\limits_{ij=12, 13, 23}\partial H_{ij}(x,\epsilon) /\partial \epsilon R_{ij}.
$$
This coordinate description might be useful, e.g., to verify (12) below. Then without loss of generality we may assume that $R$
in $so(3)$ is proportional to (say) the coordinate vectors $R_{23}$, or
$$
d_xH_{12}=d_xH_{13}=0
$$
and
$$
d_X H(x^*,0)=R_{23}\int (R(X,Y)E_2(y),E_3(y)) dy\not =0.
$$}

1) it does not equal zero,

2) it equals zero.

Below we consider in details the first (principal) case when the rank of the differential of $H(x,\epsilon)$ is at least one at
$(x^*,0)$. After that it will be easy to see that our main technical results, see the Theorems~1a and 1b below, can be obtained
in the same line of arguments also in the second case.

By the same arguments as in the fundamental lemma of the calculus of variations we see that if the rank of the map
$H(x,\epsilon)$ is one for all possible variations $\phi$ then all the vectors $R(X,Y)$ and both $\nabla_YR(X,Y)$ and $\nabla_X
R(X,Y)$ are proportional to $R$ along $\omega_{x^*}$. Indeed, take two different points $P_i=\omega_{x^*}(y_i), i=1,2$ and
assume that $R_1=I_{x^*}(y_1)R(X,Y)(P_1)$ not proportional to $R_2=I_{x^*}(y_2)R(X,Y)(P_2)$. Then choosing two $\delta$-like
functions $\phi^i_{x^*}(y)$ concentrated near these points $P_i$ we define two variations (6) for which according to (7)
derivatives of the holonomy $H(x,\epsilon)$ will be close to $R_1$ and $R_2$ correspondingly and linearly independent, making
our "taking off" possible. Hence, we come to the following conclusion.

\proclaim{Lemma~3} If for all variations $\tilde\omega_x(y)$ given by (6) the holonomy curve $H(x,\epsilon)$ goes through $id$
in $SO(3)$, then all
$$
I_{x^*}(y)R(X,Y)(y)
$$
are proportional to the vector $R$.
\endproclaim

Remark that from the Lemma~3 and the formula (7) it follows that for an arbitrary vector field $W_x(y)$ which is parallel along
$\omega_x$ we have
$$
\nabla_X W_{x^*}(y)=J^{-1}(x^*)({{L(x)}\over{2\pi}}\int\limits_{\omega_x}\psi_x(y)I^{-1}_{x^*}(y) R(X,Y)W_{x^*}(y) dy,
$$
or that
$$
\nabla_X W_{x^*}(y) \qquad \text{ is parallel to }\qquad R(X,Y)W_{x^*}(y), \tag 11
$$
where here $R(X,Y)$ is understood as the vector from $so(3)$ at $\omega_{x^*}(y)$, i.e., an anti-symmetric operator on
$\nu_{\omega_{x^*(y)}}\Sigma$.

Note, that vector fields $V(x,y)$ are parallel along $\omega_x$ for $x<x^*$ under rotations $H(x)$ which are approximately $id -
(x^*-x)R$. Because by assumption $R\not= 0$ they are close to $V^*$ such that $R(V^*)=0$, i.e., we arrive at the following
statement.

\proclaim{Lemma~4} The vector fields $V(x,y)$ parallel along $\omega_x$ and invariant under $H(x)$ converge to the vector field
$V^*$ along $\omega_{x^*}(y)$ such that
$$
R(X,Y)V^*(y)\equiv 0.
$$
The vector field $V^*(y)$ is parallel along $\omega_{x^*}$ as the limit of parallel vector fields.
\endproclaim

Next we note that taking the covariant derivative $\nabla_Y R(X,Y)$ of $R(X,Y)$ along $\omega_{x^*}$ we should have according
to the Lemma~3 the vector field which is also parallel to $R$. The same is true for the covariant derivative of the field of
transformations $R(X,Y)$ in the direction $X$ normal to $Y$. Indeed, take an arbitrary $P_1=\omega_{x^*}(y), y\not=0$ different
from $P$, and consider again a variation with a $\delta$-like function $\phi_{x^*}(y)$ concentrated near $P_1$ and zero at $P$.
Let we know that all first variations of $H(x^*,\epsilon)$ at $\epsilon=0$ are proportional to $R$. Then, as easy computations
show that the second variation of $H(x,\epsilon)$ acts on the vector $W$ as follows
$$
\delta^2_\epsilon H(x^*,\epsilon)_{|\epsilon=0}(W) =J^{-1}(x^*){{L(x)}\over{2\pi}}\int\limits_{\omega_x} \phi_x(y)I^{-1}_x(y)
(\nabla_X R(X,Y)W(y)dy +
$$
$$
J^{-1}(x^*)\int\limits_{\omega_x} \phi_x(y)I^{-1}_x(y) R(X,Y)(\nabla_X W_x(y)_{|x=x^*})dy +
$$
$$
J^{-1}(x^*)\int\limits_{\omega_x} ( {{L'(x)}\over{2\pi}} +\phi'_x(y))I^{-1}_x(y) R(X,Y)W(y)dy, \tag 12
$$
where $\phi_x'(y)$ stands for the partial derivative of $\phi_x(y)$ on $x$, $W(y)=W(x^*,y)$ and $W_x(y)$ as above is parallel
along $\omega_x$. Note that here this derivative is not a skew-adjoint map on $W$ due to the fact that $W$ also depends on the
variation, i.e., on $\epsilon$. This is given by the second right term in the previous equality which by (11) is
$$
J^{-1}(x^*)\int\limits_{\omega_x} \phi_x(y)I^{-1}_x(y) R(X,Y)(\nabla_X W_x(y)_{|x=x^*})dy \text{ proportional to } R(R(W)).
$$
Operators $R$ and $R^2$ are correspondingly the first (tangent) and the second derivatives in $GL(3)\supset SO(3)$ to the
one-parameter group of rotations in $SO(3)$ (a circle) issuing from $id$ in the same direction $R$ as the family of holonomies
$H(x)$. The third term in the right-hand side of (12) is proportional to $R(W)$. If the first vector in the right-hand side of
(12) is not proportional to $R(W)$ we again, as above in the case of the first variation non-proportional to $R$, would have
the variation which deform the holonomy curve $H(x)$ to $H(x,\epsilon)$ which for sufficiently small but positive $\epsilon$
avoids $id$ in $SO(3)$. Hence, the following is true.

\proclaim{Lemma~5} If for all variations $\tilde\omega_x(y)$ given by (6) the holonomy curve $H(x,\epsilon)$ goes through $id\in
SO(3)$, then
$$
I_{x^*}(y)\nabla_X R(X,Y)(y) \qquad \text{ and } \qquad I_{x^*}(y)\nabla_Y R(X,Y)(y)
$$
are proportional to the vector $R$.
\endproclaim

We say that the holonomy $I_{x^*}$ is vanishing along $\omega_{x^*}$ if $H(x^*)=id$ and the claims of the Lemmas~3-5 are true.
If this happens $V^*(x,y)$ belongs to the kernel of the actions of all first and second variations of $H(x^*,\epsilon)$ on
$\epsilon$ at $\epsilon=0$, the vector field $V(x,y)$ converges to $V^*(y)$ as $x\nearrow x^*$ and the following is true
$$
\nabla_X V(x,y)_{|x=x^*}\equiv 0. \tag 13
$$

Note that the same arguments work also when the vector $d_xH$ is zero. Indeed, because all $H_x$ for $x<x^*$ are non-trivial
rotations the vector field $V(x,y)$ is correctly defined, and for $x\nearrow x^*$ this vector field converges to $V^*(y)$ on
$\omega_{x^*}$. If $R(X,Y)\not\equiv 0$ then for some variation given by some profile $\phi_x(y)$ the vector $R=(\partial
H(x^*,\epsilon)/\partial \epsilon)_{|\epsilon=0}$ is not zero. Thus after a small variation given by $\phi_x(y)$ we are in the
case $R\not=0$ and may apply arguments above to obtain the claim of the Lemma~5. Similarly, if $R(X,Y)\equiv 0$ along
$\omega_{x^*}$ but the claim of the Lemma~5 is not true,  then (12) shows that after some deformation $\tilde\omega_x(y)$ we
arrive at the first case when $R\not=0$.

We summarize the obtained results in the following theorem.

\proclaim{Theorem~1a} For a domain $\Omega$ bounded by a closed curve $\omega$ there exists a smooth homotopy $\omega_x, 0\leq
x\leq 1$ of $\omega=\omega_1$ to a point $O$ such that

1) or holonomies $I_x$ along $\omega_x$ are non-trivial rotations of $\nu_{P(x)}\Sigma$ about axis $V(x)$ for all $0<x\leq 1$;

2) or for some $x\nearrow x^*<1$ these holonomies converge to the identity map. Then
$$
\nabla_X V(x,y)_{|x=x^*}\equiv 0,
$$
where $V(x,y)=I_x(y)V(x)$ are parallel translations of $V(x)$ along $\omega_x$.
\endproclaim

It might be useful to note that we may choose variations above (constructed in order to deform the initial homotopy $\omega_x$)
not diminishing domains bounded by $\omega_x$, i.e., such that we have $\Omega_x\subset \tilde \Omega_x$ for domains $\Omega_x$
$\tilde \Omega_x$ bounded by $\omega_x$ and $\tilde\omega_x$ correspondingly. Hence, varying an initial center $O$ in $\Omega$
or we can construct for a given point $O'$ a homotopy $\omega_x^{O'}$ between $\omega$ and the constant curve $O'$ with
non-vanishing holonomies along all curves in the homotopy family, or $O'\in int(\Omega^{O'})$ for some domain with vanishing
holonomy along $\partial\Omega^{O'}$.

Finely, note that for the given homotopy $\omega_x, 0\leq x\leq 1$ we may start the construction of the vector field $V(x,y)$
first by defining $V(1,y)$ and then considering the family of holonomy operators $I_x$ along $\omega_x$ for $x$ close to $1$.
In the same way as before we conclude that if $I_1$ is not the identity map then $V(x,y)$ is defined for all $x>1-\epsilon$
close enough to $1$. When the holonomy $I_1$ is trivial, i.e., $I_1=id$ but the operators $R(X,Y): V\to R(X,Y)V$ of the
infinitesimal holonomies do not vanish along $\omega_1$ we still will be able, deforming the initial homotopy $\omega_x$ as
above if necessary, to define the unit normal vector field $V(x,y)$ parallel along curves $\omega_x$ for small $x$ close enough
to $1$. Therefore, as in the Theorem~1a above we see that the following is true.

\proclaim{Theorem~1b} For a domain $\Omega$ bounded by a closed curve $\omega$ there exists a smooth homotopy $\omega_x, 0\leq
x\leq 1$ of $\omega=\omega_1$ to a point $O$ such that

1) or holonomies $I_x$ along $\omega_x$ are non-trivial rotations of $\nu_{P(x)}\Sigma$ about axis $V(x)$ for all $0<x\leq 1$;

2) or $I_1\not= id$ but for some $1>x\searrow x^*>0$ these holonomies converge to the identity map. Then
$$
\nabla_X V(x,y)_{|x=x^*}\equiv 0,
$$
where $V(x,y)=I_x(y)V(x)$ are parallel translations of $V(x)$ along $\omega_x$.

3) Or $I_1=id$ and all operators $R(X,Y): V\to R(X,Y)V$ of the infinitesimal holonomies vanish along $\omega_1$.
\endproclaim

\medskip

Below we estimate curvature of a vertical lift $\Omega_s$ of $\Omega$ in a direction of the vector field $V(x,y)$ using some
coordinates which might be different from our "coordinates" $\{x,y\}$ above.

\head \bf 2. Curvature of a local vertical lift \endhead

Consider the vertical lift of $\Omega$ along the given vector field
$$
\Omega_s(x,y)=exp_{P(x,y)}sV(x,y). \tag 14
$$
During this section the local coordinates $\{x,y\}$ in $\Omega$ will be chosen in a process of our calculations in order to
simplify them. In particular, they are not assumed to coincide with those from the previous section.

Denote by $\tilde X(x,y;s)$ and $\tilde Y(x,y;s)$ the $x$- and $y$-coordinate vectors on $\Omega_s$. By $\bar X(x,y;s)$ and
$\bar Y(x,y;s)$ we denote the horizontal lifts of $\tilde X(x,y;0)$ and $\tilde Y(x,y;0)$ (basic horizontal vector fields);
vertical $V(x,y;s)$ - the parallel transport of $V(x,y)$ along vertical geodesic from $P(x,y)$ to $P_s(x,y)=\Omega_s(x,y)$, and
by $X(x,y;s)$ and $Y(x,y;s)$ the unit vectors of the same directions as $\tilde X(x,y;s)$ and $\tilde Y(x,y;s)$. We usually
assume that at the given point of consideration (only) $\tilde X(x,y;s)$ and $\tilde Y(x,y;s)$ are unit and normal to each
other (i.e., coincide with $X(x,y;s)$ and $Y(x,y;s)$) and their first covariant derivatives vanish at this point. It holds
$$
{\Cal H}(\tilde X(x,y;s))=\bar X(x,y;s), \quad {\Cal V}(\tilde X(x,y;s))= s\nabla_X V(x,y))+o(s^2) \tag 15.1
$$
and
$${\Cal H}(\tilde Y(x,y;s))=\bar Y(x,y;s), \quad {\Cal V}(\tilde Y(x,y;s))= s\nabla_Y
V(x,y))+o(s^2). \tag 15.2
$$
Next we do calculations of some curvature tensor terms\footnote{We work first with curvature tensor terms to simplify
calculations. The obtained in Lemmas~7-9 formulas then will provide estimates for the sectional curvature of $\Omega_s$.} with
$o(s^2)$ precision, i.e., up to $O(s^2)$-terms.

\proclaim{Lemma~7}
$$
(R(\tilde X(x,y;s),\tilde Y(x,y;s))\tilde Y(x,y;s),\tilde X(x,y;s))-(R(\bar X(x,y;s),\bar Y(x,y;s))\bar Y(x,y;s),\bar X(x,y;s))
=
$$
$$ s^2((\nabla_Y(R(\bar X(x,y;0),\bar Y(x,y;0))V(x,y;0)),\nabla_X V(x,y;0))-
$$
$$
(\nabla_X(R(\bar X(x,y;0),\bar Y(x,y;0))V(x,y;0)),\nabla_Y V(x,y;0)))
$$
$$
+ s^2(R(\bar X(x,y;0),\bar Y(x,y;0))\nabla_X V(x,y;0),\nabla_Y V(x,y;0))+o(s^2).
$$
\endproclaim

\demo{Proof} We fix a point $P=P(x,y;0)$ in $\Omega$ and to simplify notations drop below $(x,y)$-arguments. By the same reason
we also drop an $s$-argument if it equals zero. By $\bar V$ we denote a vector field on $\Omega$ which is parallel at $P$,
i.e., $\nabla_X \bar V(P)=\nabla_Y \bar V(P)=0$ and equals $V$ at the point $P$. From (15) we deduce
$$
(R(\tilde X(s),\tilde Y(s))\tilde Y(s),\tilde X(s))-(R(\bar X(s),\bar Y(s))\bar Y(s),\bar X(s))
$$
$$
2s((R(\bar X(s),\bar Y(s))\bar Y(s),\nabla_X V(s))+(R(\bar X(s),\bar Y(s))\nabla_Y V(s),\bar X(s)))+
$$
$$
2s^2((R(\bar X,\bar Y)\nabla_Y V,\nabla_X V)+(R(\nabla_X V,\bar Y)\nabla_Y V,\bar X))+o(s^2)\tag 16
$$
Because
$$(R(\bar X,\bar Y)\bar Y,W)=0 \qquad \text{ and }\qquad (R(\bar X,\bar Y)U,\bar X)=0 \tag 17
$$
for any vertical $W,U$ we have
$$
(R(\bar X(s),\bar Y(s))\bar Y(s),\nabla_X V(s))+(R(\bar X(s),\bar Y(s))\nabla_Y V(s),\bar X(s))=
$$
$$
s((\nabla_V R(\bar X,\bar Y)\bar Y,\nabla_X V)- (\nabla_V R(\bar X,\bar Y)\bar X, \nabla_Y V)).
$$
By the second Bianchi identity
$$
\nabla_V R(X,Y)Y+\nabla_X R(Y,\bar V)Y+\nabla_Y R(\bar V,X)Y=0, \tag 18
$$
and (17)
$$
\nabla_V R(X,Y)Y=\nabla_Y R(X,\bar V)Y. \tag 19
$$
Or, using
$$
R(X,\bar V)Y - R(\bar V,Y)X=R(X,\bar V)Y + R(Y,\bar V)X=R(X+Y,\bar V)(X+Y) - R(X,\bar V)X - R(Y,\bar V)Y=0 \tag 20
$$
and the first Bianchi identity
$$
R(X,\bar V)Y + R(\bar V,Y)X + R(Y,X)\bar V=0 \tag 21
$$
we conclude
$$
\nabla_V R(X,Y)Y=\nabla_Y R(X,\bar V)Y=(\nabla_Y R)(X,Y)\bar V/2. \tag 22
$$
In the same way
$$
\nabla_V R(X,Y)X=\nabla_X R(Y,\bar V)X=(\nabla_X R)(X,Y)\bar V/2. \tag 23
$$
and
$$
2s((R(\bar X(s),\bar Y(s))\bar Y(s),\nabla_X V(s))+(R(\bar X(s),\bar Y(s))\nabla_Y V(s),\bar X(s)))=
$$
$$
s^2(((\nabla_Y R)(X,Y)V,\nabla_X V)-((\nabla_X R)(X,Y)V, \nabla_Y V))
$$
in (16).

Next we note that for an arbitrary operator the derivative $\nabla_X(R(V))$ equals $(\nabla_X R)(V)+R(\nabla_X V)$, or
$$
((\nabla_Y R)(X,Y)V,\nabla_X V)-((\nabla_X R)(X,Y)V, \nabla_Y V) + 2(R(X,Y)\nabla_Y V,\nabla_X V)=
$$
$$
( ( (\nabla_Y R)(X,Y)V,\nabla_X V) + (R(X,Y)\nabla_Y V,\nabla_X V)) - ( ( (\nabla_X R)(X,Y)V, \nabla_Y V) + (R(X,Y)\nabla_X
V,\nabla_Y V))
$$
$$
= (\nabla_Y (R(X,Y)V),\nabla_X V)-(\nabla_X(R(X,Y)V), \nabla_Y V). \tag 24
$$
which implies from (16) that
$$
(R(\tilde X(s),\tilde Y(s))\tilde Y(s),\tilde X(s))-(R(\bar X(s),\bar Y(s))\bar Y(s),\bar X(s))=
$$
$$
 s^2((\nabla_Y
(R(X,Y)V),\nabla_X V)-(\nabla_X(R(X,Y)V), \nabla_Y V)) + 2s^2(R(\nabla_X V,\bar Y)\nabla_Y V,\bar X)+o(s^2).\tag 16*
$$
Again, from the first Bianchi identity
$$
(R(W,Y)U,X) + (R(Y,U)W,X) + (R(U,W)Y,X)=0 \tag 25
$$
and
$$
(R(W,Y)U,X) - (R(Y,U)W,X)=(R(W,Y)U,X) + (R(U,Y)W,X)=
$$
$$
(R(W+U,Y)W+U,X)-(R(W,Y)W,X)-(R(U,Y)U,X)=0 \tag 26
$$
following from (17), we have
$$
(R(\nabla_X V,\bar Y)\nabla_Y V,\bar X))=(R(X,Y)\nabla_Y V,\nabla_X V)/2 \tag 27
$$
which finely implies through (16*)
$$
(R(\tilde X(s),\tilde Y(s))\tilde Y(s),\tilde X(s))-(R(\bar X(s),\bar Y(s))\bar Y(s),\bar X(s))=
$$
$$
 s^2((\nabla_Y(R(X,Y)V),\nabla_X V)-(\nabla_X(R(X,Y)V), \nabla_Y V)) + s^2(R(X,Y)\nabla_Y V,\nabla_X V)+o(s^2).\tag 28
$$
The Lemma~7 is proved.
\enddemo

Note that the Lemma~7 formula can be re-written as follows.

\proclaim{Lemma~8}
$$
(R(\tilde X(s),\tilde Y(s))\tilde Y(s),\tilde X(s))-(R(\bar X(s),\bar Y(s))\bar Y(s),\bar X(s)) =
$$
$$
s^2(\|R(X,Y)V\|^2-det \vmatrix
\nabla_X\nabla_X V & \nabla_X\nabla_Y V\\
\nabla_Y\nabla_X V  & \nabla_Y\nabla_Y V\\
\endvmatrix )
 + s^2 D +o(s^2), \tag 29
$$
where multiplications in a formal determinant above are scalar products and $D$ is:
$$
D= {{1}\over{2}}(- Y(\nabla_Y\nabla_XV,\nabla_XV) - X(\nabla_X\nabla_YV,\nabla_YV)+ X(\nabla_Y\nabla_XV,\nabla_YV) +
X(\nabla_Y\nabla_YV,\nabla_XV)).
$$
\endproclaim

\demo{Proof} Note that, e.g.,
$$
(\nabla_Y (R(X,Y)V),\nabla_X V)= Y(R(X,Y)V,\nabla_X V) - (R(X,Y)V,\nabla_Y \nabla_X V))
$$
Repeating similar transformations for all terms in (28) and in
$$
(R(X,Y)\nabla_Y V, \nabla_X V)=(\nabla_X \nabla_Y \nabla_Y V, \nabla_X V)-(\nabla_Y \nabla_X \nabla_Y V, \nabla_X V)
$$
we obtain for the right-hand term of (28)
$$
(\nabla_Y(R(X,Y)V),\nabla_X V)-(\nabla_X(R(X,Y)V), \nabla_Y V)) + (R(X,Y)\nabla_Y V,\nabla_X V)=
$$
$$
D-((R(X,Y)V,\nabla_Y\nabla_X V)-(R(X,Y)V,\nabla_X\nabla_Y V)+(\nabla_Y \nabla_Y V, \nabla_X \nabla_X V)-(\nabla_X \nabla_Y
V,\nabla_Y\nabla_X V))
$$
$$
D + \|R(X,Y)V\|^2 - det \vmatrix
\nabla_X\nabla_X V & \nabla_X\nabla_Y V\\
\nabla_Y\nabla_X V  & \nabla_Y\nabla_Y V\\
\endvmatrix ,
$$
since $R(X,Y)V=\nabla_X\nabla_Y V-\nabla_Y\nabla_X V$ and hence
$$
(R(X,Y)V,\nabla_Y\nabla_X V)-(R(X,Y)V,\nabla_X\nabla_Y V)=-\|R(X,Y)V\|^2;
$$
while for the "$D$-part" (which contains derivatives of scalar products) in the formula above we obtain
$$
D= {{1}\over{2}}( Y(\nabla_X\nabla_YV,\nabla_XV) - Y(\nabla_Y\nabla_XV,\nabla_XV) - X(\nabla_X\nabla_YV,\nabla_YV)+
$$
$$
X(\nabla_Y\nabla_XV,\nabla_YV) + X(\nabla_Y\nabla_YV,\nabla_XV) - Y(\nabla_X\nabla_YV,\nabla_XV) ), \tag 30
$$
which implies the claim of the Lemma~8. The Lemma~8 is proved.
\enddemo

Recall that we do calculations with $o(s^2)$-precision, i.e., up to $O(s^2)$-terms. Next we estimate the "external curvature"
term\footnote{recall the footnote before Lemma~7} of $\Omega_s$ at some point $P(x,y;s)$. To make calculations simpler we may
assume that covariant derivatives of coordinate vector fields vanish at this point, i.e.,
$$
\nabla_X X(P)=\nabla_X Y(P)=\nabla_Y Y(P)=0.
$$
Rotating if necessary the orthonormal base $\{X,Y\}$ we may also assume that vertical vectors $\nabla_X V$ and $\nabla_Y V$ are
normal at the point $P(x,y;0)$. We denote them by $dW$ and $eU$ where $\{W,U\}$ unit and normal to each other. Then the normal
space of $\Omega_s$ at the considered point $P_s(x,y)$ is generated by $\{V,\tilde M,\tilde N\}(x,y;s)$ where
$$
\tilde M=W-sd\bar X \qquad \text{ and } \qquad \tilde N=U-se\bar Y
$$
Correspondingly, the unit normals to $\Omega_s$ which are normal to each other are
$$
M=\tilde M/\|\tilde M\|={{W - s d\bar X}\over{\sqrt{1+(sd)^2}}} + o(s^2), \quad \text{ and } \quad N=\tilde N/\|\tilde N\|={{U
- s e\bar Y}\over{\sqrt{1+(se)^2}}} + o(s^2).
$$
From the Gauss equation we see that the external curvature term $R^{ext}_s(x,y)$ of $\Omega_s$ (i.e., the difference between
the curvature term $\tilde R(x,y;s)$ of the surface $\Omega_s(\epsilon)$ and the curvature tensor term $R(x,y;s)$ of the
ambient manifold $M$ in the same two-dimensional direction) equals \footnote{See last two footnotes above. To compute the Gauss
curvature we should divide these curvature terms by the area of the element $d\tilde X\wedge d\tilde Y$.}
$$
R^{ext}_s(x,y)=\sum\limits_{Z\in\{M,N\}} (\nabla_{\tilde X}{\tilde X},Z)(\nabla_{\tilde Y}{\tilde Y},Z)-(\nabla_{\tilde
X}{\tilde Y},Z)^2, \tag 31
$$
since the normal $V$ does not contribute to the Gauss formula. Because the second fundamental form of vertical fibers vanish
along $\Sigma$ for every vertical $W$ and horizontal $X$ we have
$$
\|\nabla_W X\|= O(s), \tag 32
$$
and routine calculations give
$$
\nabla_{\tilde X}{\tilde X} - H_1=\nabla_{X+sdW}{X+sdW} =sd'_xW + sd\nabla_X W + O(s^2) \tag 33
$$
where $H_1$ tangent to $\Omega_s$. In the same way
$$
\nabla_{\tilde X}{\tilde Y} - H_2=\nabla_{X+sdW}{Y+seU} =\nabla_X Y + se'_xU + se\nabla_X U + O(s^2) \tag 34
$$
and
$$
\nabla_{\tilde Y}{\tilde Y} - H_3=\nabla_{Y+seU}{Y+seU}=se'_yU + se\nabla_Y U + O(s^2) \tag 35
$$
for some $H_2, H_3$ tangent to $\Omega_s$. Which after substitution into (30) leads to the following formulas up to
$O(s^2)$-terms
$$
(\nabla_{\tilde X}{\tilde X}, M)=(sd'_xW + sd\nabla_X W,W - s d\bar X)=sd'_x,
$$
$$
(\nabla_{\tilde X}{\tilde Y}, M)=(\nabla_X Y + se'_xU + se\nabla_X U,W - s d\bar X)=(\nabla_X Y,W)+ se(\nabla_X U,W),
$$
$$
(\nabla_{\tilde Y}{\tilde Y}, M)=(se'_yU + se\nabla_Y U,W - s d\bar X)=se(\nabla_Y U,W);
$$
and
$$
(\nabla_{\tilde X}{\tilde X}, N)=(sd'_xW + sd\nabla_X W,U - s e\bar Y)=sd(\nabla_X W,U),
$$
$$
(\nabla_{\tilde X}{\tilde Y}, N)=(\nabla_X Y + se'_xU + se\nabla_X U,U - s e\bar Y)=(\nabla_X Y,U)+ se'_x,
$$
$$
(\nabla_{\tilde Y}{\tilde Y}, N)=(se'_yU + se\nabla_Y U,U - s e\bar Y)= se'_y;
$$
which with the help of
$$
\nabla_X Y= - s R(X,Y)V/2
$$
implies for the external curvature term
$$
s^{-2}R^{ext}_s(x,y)=d'_xe(\nabla_Y U,W)-(-(R(X,Y)V,W)/2+ e(\nabla_X U,W))^2 +
$$
$$
e'_yd(\nabla_X W,U)-(-(R(X,Y)V,U)/2+ e'_x)^2.
$$
Because
$$
d'_x=(\|\nabla_XV\|)'_x=(\nabla_X\nabla_X V,\nabla_X V)/\|\nabla_X V\|=(\nabla_X\nabla_X V,W)
$$
$$
e'_x=(\|\nabla_YV\|)'_x=(\nabla_X\nabla_Y V,\nabla_Y V)/\|\nabla_Y V\|=(\nabla_X\nabla_Y V,U)
$$
$$
e'_y=(\|\nabla_XV\|)'_y=(\nabla_Y\nabla_Y V,\nabla_Y V)/\|\nabla_Y V\|=(\nabla_Y\nabla_Y V,U)
$$
and, e.g.,
$$
(\nabla_Y U,W)=\|\nabla_Y V\|(\nabla_Y {{\nabla_Y V}\over{\|\nabla_Y V\|}},W)=(\nabla_Y \nabla_Y V,W)
$$
due to the fact that $(W,U)=0$; by direct calculations we conclude
$$
s^{-2}R^{ext}_s(x,y)=
$$
$$
(\nabla_X\nabla_X V,W)(\nabla_Y\nabla_Y V,W)-((R(X,Y)V,W)^2/4+(R(X,Y)V,W)(\nabla_X\nabla_Y V,W)-(\nabla_X\nabla_Y V,W)^2)+
$$
$$
(\nabla_X \nabla_X V,U)(\nabla_Y \nabla_Y V,U) - (R(X,Y)V,U)^2/4 +(R(X,Y)V,U)(\nabla_X \nabla_Y V,U) - (\nabla_X \nabla_Y
V,U)^2)=
$$
$$
(\nabla_X\nabla_X V,\nabla_Y\nabla_Y V)-\|R(X,Y)V\|^2/4+(R(X,Y)V,\nabla_X\nabla_Y V) - \|\nabla_X\nabla_Y V\|^2. \tag 36
$$
Because
$$
(R(X,Y)V,\nabla_X\nabla_Y V)=(\nabla_X \nabla_Y V-\nabla_Y\nabla_X V,\nabla_X\nabla_Y V)
$$
we obtain the following statement.

\proclaim{Lemma~9} The external curvature term $R^{ext}_s(x,y)$ of $\Omega_s$ is given by the formula
$$
R^{ext}_s(x,y)=s^2(-\|R(X,Y)V\|^2/4 + det \vmatrix
\nabla_X\nabla_X V & \nabla_X\nabla_Y V\\
\nabla_Y\nabla_X V  & \nabla_Y\nabla_Y V\\
\endvmatrix )
$$
where, as before, multiplications in a formal determinant above are scalar products.\footnote{We may complete the proof of the
Lemma~9 with the following analysis of the external curvature term of $\Omega_s$ at the point $O(s)=P(x,y;0)$ in the
coordinates from the Theorem~1a, where our $y$-"coordinate" curves $\omega_x$ degenerate to a point and $kg(x,y)\to\infty$ as
$x\to 0$. Easy to see this is singularity of the coordinate system, which does not yield the singularity of $\Omega_s$. Indeed,
as we already noted, at this point the vector field $V$ is smooth with all first-order covariant derivatives vanishing, i.e.,
$d(0,y)=0$. Take other than our "polar"-type coordinates: let, for instance, the "new" coordinates $\{x',y'\}$ be the normal
coordinates on $\Sigma$ with the center at $O$. Because $\omega_x$ for small $x$ are circles with radius $x$ around $O$ we have
in new coordinates:
$$
\nabla_{Y cos(\alpha) - X sin(\alpha)} V(x cos(\alpha), x sin(\alpha))\equiv 0
$$
for all $0\leq \alpha\leq 2\pi$, which after taking derivative on $\alpha$ gives
$$
\nabla_{X}\nabla_{X}V(0,0), \nabla_{X}\nabla_{Y}V(0,0), \nabla_{Y}\nabla_{Y}V(0,0)=0,
$$
which in turn as in (15) implies
$$
Y'(0,0;s)=\bar Y(0,0;s), \quad X'(0,0;s))=\bar X(0,0;s) \qquad \text{ and } \qquad \nabla_{X'}{X'}(x',y';s),
\nabla_{Y'}{Y'}(x,y;s) = 0 \tag 15'
$$
for new coordinate vectors $X'(x',y';s), Y'(x',y';s)$ on $\Omega_s$, and by (31) implies for the external curvature of
$\Omega_s$ at $P(0,0;s)$ the claim of the Lemma. The proof of the Lemma~9 is complete.}
\endproclaim

Now we put formulas above together and draw some conclusions. Because by the fundamental O'Neill's formula (4)
$$
(R(\bar X(s),\bar Y(s))\bar Y(s),\bar X(s))=(R(\bar X,\bar Y)\bar Y,\bar X) - 3\|A_XY\|^2=
$$
$$
(R(\bar X,\bar Y)\bar Y,\bar X) - {{3}\over{4}}\|R(X,Y)V\|^2
$$
and from the Gauss fundamental equation for the curvature term $\tilde R(x,y;s)$ of the surface $\Omega_s$ follows
$$
\tilde R(x,y;s)=(R(\tilde X(x,y;s),\tilde Y(x,y;s))\tilde Y(x,y;s),\tilde X(x,y;s)) + R^{ext}_s(x,y),
$$
from Lemmas~7-9 we conclude our main formula of this section
$$
\tilde R(x,y;s) = \tilde R(x,y;0) + s^2 D + o(s^2), \tag 37
$$
where
$$
D= {{1}\over{2}}(- Y(\nabla_Y\nabla_XV,\nabla_XV) - X(\nabla_X\nabla_YV,\nabla_YV)+ X(\nabla_Y\nabla_XV,\nabla_YV) +
X(\nabla_Y\nabla_YV,\nabla_XV)), \tag 38
$$
and the curvature term $\tilde R(x,y;0)$ of $\Omega_0$ equals $(R(X,Y)Y,X)(x,y;0)$.\footnote{Similar but simpler calculations
give
$$
{((R(\tilde X,\tilde Y)\bar Y,\bar X)(x,y;s))'_s}_{|s=0}={((R(\tilde X,\tilde Y)\bar Y,\bar X)(x,y;s))"_s}_{|s=0}=0.
$$}
The sectional curvature $\tilde K(x,y;s)$ of $\Omega_s$ at the point $P(x,y;s)$ is given through the curvature term as
$$
\tilde K(x,y;s)={{\tilde (R(\tilde X(x,y;s),\tilde Y(x,y;s))\tilde Y(x,y;s),\tilde X(x,y;s))}\over{\|\tilde X(x,y;s)\wedge
\tilde Y(x,y;s)\|^2}}.
$$
Because
$$
\|\tilde X(x,y;s)\wedge \tilde Y(x,y;s)\|^2=1 + s^2(d^2+e^2) \tag 39
$$
the equality (37) gives for the curvature forms
$$
\tilde K(x,y;s) d\sigma_s = \tilde K(x,y;0)(1-{{s^2}\over{2}}(d^2+e^2)) d\sigma_0 + s^2 Dd\sigma_0 + o(s^2), \tag 40
$$
where $d\sigma_s$ denotes the area form $d\tilde X(x,y;s)\wedge d\tilde Y(x,y;s)$ of $\Omega_s$.

\head \bf 3. Local parallel section. Vanishing holonomy case \endhead

At this point we should repeat again that to simplify Lemmas~7-9 calculations we have used local coordinates in $\Omega$
satisfying some assumptions such as: coordinate vectors $X(x,y)$ and $Y(x,y)$ at the given point $P$ were orthonormal and such
that derivatives of the vector field $V(x,y)$ given by the Theorem~1 $\nabla_X V$ and $\nabla_Y V$ were orthogonal. We also
assumed above that
$$
\nabla_X X(P)=\nabla_X Y(P)=\nabla_Y Y(P)=0 \tag 41
$$
at the point $P$ where we calculated terms of (40). However in the obtained formula (40) not only curvature terms do not depend
on this particular choice of coordinates, but also the term
$$
d^2+e^2=(\nabla_X V,\nabla_X V)+(\nabla_Y V,\nabla_Y V) \tag 42
$$
can be rewritten as an invariant, known as the "vertical" part of the energy of our vertical lift $V:P(x,y)\to V(x,y)$, as
follows:
$$
E^{\Cal V}(V(x,y))=g^{ij}(x,y)(\nabla_i V,\nabla_j V) \tag 43
$$
for an arbitrary coordinate system $\{x^1,x^2\}$ in the neighborhood of $P$ in $\Sigma$, where $\nabla_i
V=\nabla_{\partial/\partial x^i}V$ and $g_{ij}(x^i,x^j)$ and $g^{ij}(x^i,x^j)$ are metric tensor and its inverse
correspondingly. Hence, the same is true also for the $D$-term in (40): it can be expressed in a form which is invariant under
coordinate changes. The exact formula easily follows from its origin from Lemmas~7-9's calculations and is left until the next
paper where we study its properties in more details. In this paper the following property of $D$ is crucial.

\proclaim{Lemma~10} The two-form $D(x,y)dX(x,y)\wedge dY(x,y)$ is exact
$$
D(x,y)dX(x,y)\wedge dY(x,y)=d\eta(x,y),
$$
where the one-form $\eta$ has the type:
$$
\eta(x,y)= A(x,y)dX + B(x,y)dY,
$$
with coefficients $A(x,y), B(x,y)$ of the form
$$
A(x,y)=(A_1(x,y), \nabla_{X} V(x,y))+(A_2(x,y), \nabla_{Y} V(x,y))
$$
and
$$
B(x,y)=(B_1(x,y), \nabla_{X} V(x,y))+(B_2(x,y), \nabla_{Y} V(x,y))
$$
\endproclaim

\demo{Proof} The proof is immediate by the definition of the differential:
$$
D(x,y)dX(x,y)\wedge dY(x,y) = (Y(\nabla_X\nabla_YV,\nabla_XV) - Y(\nabla_Y\nabla_XV,\nabla_XV) -
X(\nabla_X\nabla_YV,\nabla_YV)+
$$
$$
X(\nabla_Y\nabla_XV,\nabla_YV) + X(\nabla_Y\nabla_YV,\nabla_XV) - Y(\nabla_X\nabla_YV,\nabla_XV)) dX(x,y)\wedge dY(x,y) =
$$
$$
d((\nabla_X\nabla_YV,\nabla_XV) dX) - d((\nabla_X\nabla_YV,\nabla_yV)dY) + d((\nabla_Y\nabla_XV,\nabla_YV)dY) +
d((\nabla_Y\nabla_YV,\nabla_XV)dY) =
$$
$$
{{1}\over{2}} d ((d^2)'_ydX - (e^2)'_xdY) - d( (f^2)'_xdX - (f^2)'_ydY ),
$$
where  $f^2=(\nabla_XV,\nabla_YV)$; or
$$
\eta ={{1}\over{2}} ((d^2)'_ydX - (e^2)'_xdY) - ((f^2)'_xdX - (f^2)'_ydY ). \tag 44
$$
\enddemo

\medskip

Now we can prove our main technical results. The first one is about the vector field $V(x,y)$ from the Theorem~1a. Denote for
short by $\Omega^*=\{P(x,y) | x\leq x^*\}\subset \Omega$ the domain where the vector field $V(x,y)$ is defined and by
$\Omega^*_s$ the vertical lift of $\Omega^*$ in direction of this vector field.

\proclaim{Theorem~2a} When the holonomy vanishes along the boundary $\omega_{x^*}$ of $\Omega^*$ then the vector field $V(x,y)$
constructed in the Theorem~1a is parallel on $\Omega^*$:
$$
\nabla_X V(x,y)=\nabla_Y V(x,y) \equiv 0.
$$
\endproclaim

\demo{Proof} First we note that all the curves $\omega^*_{s}(y)=\partial\Omega^*_s$ which are vertical lifts of
$\omega_{x^*}(y)$ have the same geodesic curvature. Indeed, from (1) it follows $\nabla_V Y, [V,Y]\equiv 0$ and $R(V,Y)Y\equiv
0$. Hence, from
$$
\nabla_V\nabla_Y Y = R(V,Y) + \nabla_Y\nabla_V Y + \nabla_{[V,Y]} Y = R(V,Y)Y
$$
we have
$$
\nabla_Y Y (\omega^*_s(y)) \equiv kg(y) \bar X(\omega^*_s(y)), \tag 45
$$
where $kg(y)$ stands for the geodesic curvature of $\omega_{x^*}(y)$. In the case of the vanishing holonomy by the Theorem~1
and (15) the tangent subspace to $\Omega^*_s$ along $\omega^*_s(y)$ coincides with the horizontal subspace, i.e., contains the
vector $\nabla_Y Y (\omega^*_s(y))$ of the geodesic curvature of this vertical lift of $\omega_{x^*}$, which implies that the
geodesic curvature of $\omega^*_s(y)$ in $\Omega^*_s$ is the same as the geodesic curvature of $\omega_{x^*}(y)$ in $\Omega^*$.
Hence, by the Gauss-Bonnet theorem
$$
\int\limits_{\Omega^*_s}\tilde K(x,y;s) d\sigma_s =\int\limits_{\Omega^*_s} \tilde K(x,y;0)d\sigma_0. \tag 46
$$
If we compare this with (40) we get
$$
\int\limits_{\Omega^*} \tilde K(x,y;0)(d^2+e^2) d\sigma_0= 2 \int\limits_{\Omega^*} D d\sigma_0. \tag 47
$$
By the Lemma~10 and Stokes theorem
$$
\int\limits_{\Omega^*} D d\sigma_0=\int\limits_{\omega_{x^*}} \eta, \tag 48
$$
which in turn equals zero since by the Theorem~1 the one-form $\eta$ vanishes identically along $\omega_{x^*}$. I.e., we have
$$
\int\limits_{\Omega^*} \tilde K(x,y;0)(d^2+e^2) d\sigma_0=0,
$$
or
$$
\tilde K(x,y;0)(d^2+e^2)\equiv 0 \tag 49
$$
from the non-negativity of the curvature. By the "prism"-construction the holonomy $I_c$ along a closed curve $c$ vanishes if
$c$ is inside some open domain in $\Sigma$ with zero curvature and is contractible in this domain, see the Lemma~3.6 [M3].
Therefore, $d$ and $e$ vanish in the interior of the closure of the set in $\Omega$ where $K(x,y;0)$ equals zero. Which leads
to
$$
d(x,y)=e(x,y)=0 \qquad \text{ if }\qquad K(x,y;0)=0 \tag 50
$$
because $d(x,y), e(x,y)$ are smooth functions, and
$$
d(x,y)=e(x,y)=0 \qquad \text{ for all }\qquad P(x,y)\in\Omega^* \tag 51
$$
with the help of (50). The Theorem~2a is proved.
\enddemo

Note that the proven result does not mean that the holonomy on $\Omega^*$ is trivial. We have proved only that the vector field
$V(x,y)$ constructed in the Theorem~1a is parallel on $\Omega^*$, which does not imply that the infinitesimal holonomy
operators $R(X,Y)(x,y)$ vanishes identically. Note also that under condition: $R(X,Y)(x,y)\not=0$ the vector field $V(x,y)$
coincides with $V^*(x,y)$ (which is not defined otherwise).

\medskip

Next we note that the form $\eta$ also vanishes along an arbitrary geodesic: if some curve $\omega_x(y)$ is a geodesic and
$V(x,y)$ is a vector field which is parallel along $\omega_x$ then in a local half-geodesic coordinate system with $\omega_x$
as an axe (such system of coordinates obviously satisfies our restrictions on coordinate systems where the form $\eta$ is given
by the formula (44) above) it holds:
$$
\eta(\dot\omega_x(y))= - \eta(\partial/\partial y)=-(e^2)'_x=(\nabla_X\nabla_Y V(x,y) Y,\nabla_YV(x,y))=0 \tag 52
$$
since $\nabla_Y V(x,y)=0$ by the definition of $V(x,y)$. This implies our second main technical result.

\proclaim{Theorem~2b} If $\Omega$ is bounded by the closed geodesic $\omega(y)$, and we have 1) or 2) case in the Theorem~1b
then the vector field $V(x,y)$, which existence is stated in the Theorem~1b in $\Omega$ or $\Omega\backslash \Omega_{x^*}$, is
parallel in the corresponding domain.
\endproclaim

\demo{Proof} The proof is immediate by the same arguments as in the proof of the Theorem~2a. If the vector field $V(x,y)$ is
defined on the whole $\Omega$ we can define the vertical lift $\Omega_s$. Because the boundary $\omega$ of $\Omega$ is a
geodesic its vertical lifts $\omega_s$ are also geodesics in $M$ by (45). Then as in the Theorem~2a by the Gauss-Bonnet theorem
it holds (46) which with the help of the Stokes formula and (52) implies (49) and the claim (51) of our theorem as above. If
the vector field $V(x,y)$ is defined only on some sub-domain $\tilde\Omega(x^*)=\Omega\backslash \Omega_{x^*}$ for $0<x^*<1$
(i.e., we have the second case in the Theorem~1b), then we apply (45) and the Gauss-Bonnet theorem to the vertical lift
$\tilde\Omega_{s}(x^*)$ of this sub-domain. Then the Lemma~10 together with (52) infer (51).
\enddemo

\head \bf 4. Global section. The case of non-vanishing holonomy \endhead

If the holonomy never vanishes we may, actually, construct a global parallel section $V:\Sigma\to \nu\Sigma$ of the unit normal
bundle of $\Sigma$. The proof is easy by going to contradiction. Indeed, assume that at some point $O_+\in\Sigma$ the holonomy
operator $R(X,Y)$ is not zero. Then, as we have seen already, in the neighborhood of this point $O_+$ the smooth vector field
$V^*$ is correctly defined. Assume that it is not parallel, i.e.,
$$
\nabla_X V^*\qquad \text{ or }\qquad \nabla_X V^* \qquad\text{ not zero.} \qquad \tag 53
$$
Take another point $O_-$, a disk $\Omega^r$ with a center $O_-$ of a small radius $r$. Next consider the homotopy $\omega_x$ of
the boundary of this disk $\omega^r$ to a point $O_+$ inside $\Omega=\Sigma\backslash \Omega^r$. Then the family of parallel
transports $I_x$ along $\omega_x$ never vanishes for otherwise we would not have (52) by the Theorem~2. Thus, taking $r\to 0$
we can define the vector field $V$ as in the Theorem~1 on $\Sigma$.\footnote{First on $\Sigma\backslash\{O_-\}$, but then
arguments as above in footnote~10 shows that the vector field $V$ can be smoothly continued to the point $O_-$ as well.} Let us
call this $(O_+,O_-)$-homotopy.

Now, applying Lemmas~7-9 computations to this global section $V$ instead of (48) we have
$$
\int\limits_{\Sigma} D d\sigma_0=\int\limits_{\partial\Sigma} \eta=\int\limits_{\emptyset} \eta=0, \tag 54
$$
implying (51) as before, i.e., that the constructed section is parallel. Therefore, the following is true.

\proclaim{Lemma~11} If the holonomy does not vanish for any $(O_+,O_-)$-homotopy then there exists a global parallel section
$V:\Sigma\to \nu\Sigma$ of a unit normal bundle such that the family of corresponding lifts
$$
\Sigma_s=\{exp_P sV(P) \quad | \quad P \in\Sigma, \quad | \quad 0\leq s \leq s_0\}
$$
is isometric to the direct product $\Sigma\times [0,s_0]$.
\endproclaim

When the global parallel section $V:P\in\Sigma \to V(P)\in\nu\Sigma$ exists the proof of the Theorem~A is easy. Indeed, then
all horizontal lifts $\Sigma_s$ are totally geodesic sub-manifolds in $M$ isometric copies of $\Sigma$, or pseudo-souls. Thus,
arguing in the same way as in the original paper by Cheeger and Gromoll, see [CG] or [Y], we can prove that the sectional
curvature of $M$ vanish in all two-dimensional "vertizontal" directions along $\Sigma_s$, i.e., generated by one vector tangent
to $\Sigma_s$ and another - normal to it.\footnote{or even, using Perelman's arguments, that (1) is fulfilled for $\Sigma_s$
too. For the proof consider "up-and-down" construction from [M1] and proceed as in [P].}

It would be interesting to understand when the global section exists. Note that, as we will prove in an instant (see the next
section) along every closed geodesic $\gamma$ on $\Sigma$ their exists a  parallel normal vector field with vanishing covariant
derivatives. Thus, it would be natural to conjecture the existence of the global parallel section in a case when through every
point of $\Sigma$ goes some closed geodesic.

\head \bf 5. The proof of the Theorem~A \endhead

There exists at least one closed geodesic $\gamma$ in $\Sigma$ which is contractible since $\Sigma$ is diffeomorphic to the
sphere $S^2$. Consider the homotopy $\omega_x, 0\leq x\leq 1$ between some point $O$ and $\gamma$. According to the Theorems~1b
and 2b there exists a vector field $V^*(y)$ parallel along $\gamma(y)$ such that
$$
R(X,Y)V^*(y)\equiv 0. \tag 55
$$
Then by the Lemmas~1 and 2 we conclude
$$
A(Q(V^*(y),s))\equiv 0, \tag 56
$$
i.e., along all geodesics $\gamma_s(y)=exp_{\gamma(y)}sV^*(y)$ which are horizontal lifts of $\gamma$ the fundamental
$A$-tensor vanishes. This implies that the family of vertical spaces ${\Cal V}(Q(V^*(y),s))$ is parallel along $\gamma_s$. Fixe
some $s>0$, two unit and vertical parallel vector fields $W(y)$ and $U(y)$ along $\gamma_s(y), 0\leq y\leq 2\pi$, and consider
the mean curvature vector $H(y)$ of the vertical fiber at $\gamma_s(y)$:
$$
H(y)=T_{W(y)} W(y) + T_{U(y)} U(y).
$$
It does not depend on the particular choice the orthonormal base $\{W(y),U(y)\}$ and therefore is a smooth vector field along
$\gamma_s$. Thus the scalar product $(H(y),Y(y))$, where $Y(y)=\dot\gamma_s(y)$ is a periodic function along $\gamma_s$. For
its derivative we have
$$
(H(y),\dot\gamma_s(y))'_y=( (\nabla_YT)_{W(y)} W(y) + (\nabla_YT)_{U(y)} U(y), Y(y) )  \tag 57
$$
since $\gamma_s$ is a geodesic and $\{W(y),U(y)\}$ are parallel along $\gamma_s$. An integral of (57) over closed $\gamma_s$
equals zero, which with the help of (56) and (3) implies
$$
\int\limits_{\gamma_s}(R[W(y),Y(y)] + R[U(y),Y(y)]) dy = - \int\limits_{\gamma_s}(\|T_W Y(y)\|^2+\|T_U Y(y)\|^2 ) dy. \tag 58
$$
Because the curvature is non-negative we conclude that curvatures $R[W(y),Y(y)]$ and $R[U(y),Y(y)]$ vanish along the geodesic
$\gamma_s$ together with the second fundamental form of vertical fibers relative to the normal $Y(y)$:
$$
R[W(y),Y(y)]=R[W(y),Y(y)]\equiv 0 \tag 59
$$
and
$$
T_W(y) Y(y)=T_U(y) Y(y)\equiv 0, \tag 60
$$
which due to the Gauss fundamental equation implies that not only the sectional curvature of $M$ vanish along $\gamma_s$ in
two-dimensional directions $\{W(y),Y(y)\}$ but also that the sectional curvature of the hypersurface $N\Sigma$ in the same
direction equals zero.

Theorem~A is proved.

\enddocument

\Refs
\widestnumber \key {AAAA}


\ref \key CG \by J. Cheeger, D. Gromoll \paper On  the  structure  of  complete manifolds of nonnegative curvature \jour Ann.
Math. \vol 96 no.3\yr 1972 \pages 413--443 \endref

\ref \key GT \by D.~Gromoll and K.~Tapp \paper Nonnegatively curved metrics on $S^2 \times R^3$\jour Geometriae Dedicata \yr
2003 \vol 99 \pages 127--136 \endref

\ref \key M1 \by V.~Marenich \paper Metric structure of open manifolds of nonnegative curvature \jour Doklady Acad. Sc. USSR
\vol 261:4 \yr 1981 \pages 801-804 \endref

\ref \key M1 rus \by V.~Marenich \paper Metric structure of open manifolds of nonnegative curvature (complete version in
russian) \jour Ukrainian Geom. Sb. \vol 26 \yr 1983 \pages 79-96 \endref

\ref \key M2 \by V.~Marenich \paper The metric of nonnegative curvature on the tangent bundle of two-dimensional sphere, \jour
Sibirsk. Math. Zh. \vol 27:2 \yr 1986 \pages 121-138 \endref

\ref \key M3 \by V.~Marenich \paper The holonomy in open manifolds of nonnegative curvature \yr 1993 \jour MSRI, Preprint
No.~003--94. \endref

\ref \key M4 \by V.~Marenich \paper The holonomy in open manifolds of nonnegative curvature \jour Michigan Math. Journal \vol
43:2 \yr 1996 \pages 263--272 \endref

\ref \key O'N1 \by B.~O'Neill \paper The fundamental equations of submersion \jour Mich. Math. J. \vol 13 no. 4 \yr 1966 \pages
459--469 \endref

\ref \key P \by G.~Perelman \paper Proof of the soul conjecture of Cheeger and Gromoll \jour J. Differential Geometry \vol 40
\yr 1994 \pages 209--212 \endref

\ref \key Y \by J.~W.~Yim \paper Spaces of souls in a complete open manifold of nonnegative curvature \jour J. Differential
Geometry \vol 32 no. 2\pages 429--455 \yr 1990 \endref

\endRefs
\bye